\documentclass[12pt, 
a4paper
]{article}
\usepackage[utf8]{inputenc}
\usepackage{amsfonts,amssymb,amsthm,amsmath}

\theoremstyle{plain}
\newtheorem*{maintheorem}{Theorem}

\theoremstyle{definition}

\theoremstyle{remark}

\newcommand{\N}{\mathbb{N}}
\newcommand{\Z}{\mathbb{Z}}
\newcommand{\R}{\mathbb{R}}

\newcommand{\de}{\mathrm{d}}

\newcommand{\id}{\operatorname{id}}

\newcommand{\interior}[1]{\mathring{#1}}

\newcommand{\fuermich}[1]{}
\newcommand{\ifindoubt}[1]{}
\newcommand{\suppressed}[1]{}

\fuermich{}

\usepackage{enumitem}

\title{Removing parametrized rays symplectically}
\author{B. Stratmann}
\date{}
\begin{document}
\maketitle
\begin{abstract}
Let $(M,\omega)$ be a symplectic manifold.
Let $[0,\infty)\times Q\subset\R\times Q$ be considered as parametrized rays $[0,\infty)$ and let $\varphi:[-1,\infty)\times Q\to M$ be an injective, proper, continuous map immersive on $(-1,\infty)\times Q$.
If for the standard
vector field $\frac{\partial}{\partial t}$ on $\R$ and any further vector field
$\nu$ tangent to $(-1,\infty)\times Q$ the equation $\varphi^*\omega(\frac{\partial}{\partial t},\nu)=0$ holds then $M$ and $M\setminus \varphi([0,\infty)\times Q)$ are symplectomorphic.
\end{abstract}
The question which subsets $N$ of a symplectic manifold $M$ can be chosen
such that $M$ and $M\backslash N$ are symplectomorphic
has been treated in particular for $M=\R^{2n}$ a time ago already, see e.g. \cite{Gro85,McD87,McDTra93,Tra93}. More recently, X.~Tang showed that for a general
manifold $M$ the subset $N$ can be chosen to be a ray if the ray possesses a ``wide neighborhood'' (\cite{Tan18}).
Roughly speaking a ray is a 2-ended connected non-compact 1-dimensional local submanifold whose one end closes up inside $M$ while at the other end the embedding is proper. 
In this paper
an extension to higher dimensional sets regarded as parametrized rays is provided. 
While 
for those higher dimensional sets a condition is needed, this condition is trivially fulfilled for
an isolated ray
as treated in \cite{Tan18}.\\
In order to state the theorem precisely
let $\frac{\partial}{\partial t}$ denote the standard vector field on $\R$, i.e. whose flow consists of translations.
\begin{maintheorem}
Let $(M, \omega)$ be a symplectic manifold, $Q$ some manifold and the map
$\varphi:[-1,\infty)\times Q\to M$ be injective, proper and continuous such that
$\varphi\vert_{(-1,\infty)\times Q}$ is immersive. If the equation
\begin{equation}
 \imath_{\frac{\partial}{\partial t}}\varphi^*\omega=0\label{condition_ContractionWithOmegaVanish}
\end{equation}
holds on $(-1,\infty)\times Q$ then $M$ and $M_0=M\backslash \varphi([0,\infty)\times Q)$ are symplectomorphic.\\
Additionally this symplectomorphism
can be chosen to be the identity
outside some selected neighborhood
of $\varphi([-1,\infty)\times Q)$.
\end{maintheorem}
The idea of the proof is to provide a
diffeotopy of symplectomorphisms 
$\psi_s:M\to M$ arising from a time-dependent vector field $\xi_s$ such that
\begin{enumerate}[label=($\xi$\arabic*)]
\item the set
\begin{equation*}
\bigcup_{s'\in [s-1,s+1]}\{m\in M\ \vert \ \xi_{s'}\vert_m\not=0\}
\end{equation*}
is relatively compact for all $s\in\R$,
\label{condition_psi_integrability}
\item for all compact subsets $K\subset M_0$ there is $s_K\in \R$ such that
$\psi_s\vert_K = \psi_{s'}\vert_K$
for all $s, s'\geq s_K$ and
\label{condition_psi_s_locally_stable}
\item for each compact set $L\subset M$
there exists $\sigma_L\in\R$ such that
\begin{equation*}
\psi_\sigma\vert_{\psi_{\sigma_L}^{-1}(L)} =\psi_{\sigma'}\vert_{\psi_{\sigma_L}^{-1}(L)}
\end{equation*}
for all $\sigma,\sigma'\geq \sigma_L$.
\label{property_psi_s_becoming_stable_in_image}
\end{enumerate}
This curve of symplectomorphisms will arise from a time-dependent Hamiltonian
vector field $\xi_s$ which in turn will be the extension of a time-dependent vector field $\zeta_s$ on $\varphi((-1,\infty)\times Q)$ arising from a diffeotopy. 
The condition $\imath_{\frac{\partial}{\partial t}}\varphi^*\omega=0$
in the Theorem (Condition~\eqref{condition_ContractionWithOmegaVanish})
imposes few constraints on the rank of $\varphi^*\omega$ such that there is
little hope to extend $\zeta_s$ to $\xi_s$ using local neighborhood models as classically known to exist around isotropic submanifolds. Therefore
integrabilty of the Hamiltonian flow
is requested explicitly by Property~\ref{condition_psi_integrability}.
Further Property~\ref{condition_psi_s_locally_stable} will ensure that the curve of symplectomorphism $\psi_s$ becomes locally stable on $M_0$ for $s\to\infty$ and therefore converges to a limit $\psi:M_0\to M$ still satisfying $\psi^*\omega=\omega$ while
Property~\ref{property_psi_s_becoming_stable_in_image} will provide surjectivity of this limit $\psi$.
\begin{proof}
In the first part of the proof a suitable diffeotopy
\begin{equation*}
\theta_s:[-1,\infty)\times Q\to [-1,\infty)\times Q 
\end{equation*}
is constructed which 
in turn is 
constructed from a diffeotopy
$\tau_s$ of $[-1,\infty)$.\\
Fix $b\in (-1,0)$ and a diffeomorphism
$\tau:[-1,0)\to[-1,\infty)$ with
$\tau\vert_{[-1,b]}=\id_{[-1,b]}$. There is a diffeotopy $\tau_s$, i.e. a smooth curve $s\mapsto \tau_s$ of
diffeomorphisms $\tau_s:[-1,\infty)\to [-1,\infty)$, such that
\begin{enumerate}[label=($\tau$\arabic*)]
\item $\tau_s=\id$ for all $s\leq 0$,
\item $\tau_s\vert_{[-1,b]}=\id_{[-1,b]}$ for all $s\in\R$,
\item for all $s\in\R$ the set
\begin{equation*}
\bigcup_{s'\in[s-1,s+1]}\{t\in(-1,\infty)\vert \tau_{s'}(t)\not =t\}
\end{equation*}
is relatively compact in $(-1,\infty)$ and
\item for each compact subset $A\subset
(-1,0)$ there is $s_A\in \R$ such that
$\tau_s\vert_A=\tau\vert_A$
for all $s\geq s_A$.
\end{enumerate}
In order to define a suitable diffeotopy
$\theta^\circ_s:[-1,\infty)\times Q\to [-1,\infty)\times Q$
choose a function $\rho:Q \to (0,\infty)$ such that the set
$\{q\in Q\ \vert\ \rho(q)\leq c\}$
is compact for each $c\in\R$.
Define $\theta^\circ_s$ by
\begin{equation*}
\theta^\circ_s(t,q)=(\tau_{s-\rho(q)}(t), q)\quad \text{as well as}\quad \theta(t,q)=(\tau(t),q)\ .  
\end{equation*}
Denoting $\pi_Q:\R\times Q\to Q$ the projection the maps $\theta^\circ_s$
satisfy $\pi_Q\circ\theta^\circ_s=\pi_Q$
while for fixed $q$ and $s$ the map $t\mapsto \tau_{s-\rho(q)}(t)$ is a diffeomorphism of $[-1,\infty)$. Thus $\theta^\circ_s$ is a diffeomorphism. By construction it satisfies
\begin{enumerate}[label=($\theta$\arabic*)]
\item $\theta^\circ_s=\id$ for all $s\leq 0$,\label{property_theta_first}
\item $\theta^\circ_s\vert_{[-1,b]\times Q}=\id_{[-1,b]\times Q}$ for all $s\in\R$,
\item for all $s\in\R$ the set
\begin{equation*}
\bigcup_{s'\in [s-1,s+1]}\{(t,q)\in (-1,\infty)\times Q\ \vert\ \theta^\circ_{s'}(t,q)\not=(t,q)\}
\end{equation*}
is relatively compact in $(-1,\infty)\times Q$ and
\item for each compact set $A\subset [-1,0)\times Q$ there is $s_A\in\R$ such that $\theta^\circ_s\vert_A=\theta\vert_A$ for all $s\geq s_A$.\label{property_theta_last}
\label{property_theta_becoming_stable}
\end{enumerate}
Next, the time $s$ of the diffeotopy 
will be deformed
which will be helpful later
to cut off time-dependent functions (see~\ref{condition_chi_cutoffArea}). Choose an increasing
smooth map $\kappa:\R\to\R$ and
$\delta>0$ such that
$\kappa\vert_{[n-\delta,n+\delta]}=n$ for all $n\in\Z$
and define the diffeotopy $\theta_s=
\theta^\circ_{\kappa(s)}$
satisfying likewise all the above
properties~\ref{property_theta_first}-\ref{property_theta_last}.
This diffeotopy can be seen to be the flow of a time-dependent vector field $\zeta_s$ on $[-1,\infty)\times Q$. The property $\pi_Q\circ\theta_s=\pi_Q$ shows that the time-dependent vector field $\zeta_s$ 
points in the direction of the rays, i.e. there is a time-dependent
function $\lambda_s:[-1,\infty)\times Q\to \R$ such that 
\begin{equation}
\zeta_s=\lambda_s\cdot\frac{\partial}{\partial t}\ .\label{equation_pointingInDirectionOfRays}
\end{equation}
Furthermore $\zeta_s$ satisfies the following properties.
\begin{enumerate}[label=($\zeta$\arabic*)]
\item $\zeta_s=0$ for all $s\leq 0$ and for all $s\in [n-\delta,n+\delta]$
for each $n\in\Z$,
\label{property_zeta_constantAroundN}
\item $\zeta_s\vert_{[-1,b]\times Q}=0$
for all $s\in\R$,
\item for each $s\in\R$ the set
\begin{equation*}
\bigcup_{s'\in[s-1,s+1]}\{(t,q)\in (-1,\infty)\times Q\ \vert\ \zeta_{s'}\vert_{(t,q)}\not=0\}
\end{equation*}
is relatively compact in $(-1,\infty)\times Q$ and\label{property_zeta_integrabilty}
\item for each compact subset
$A\subset [-1,0)\times Q$ there is an
$s_A\in\R$ defining $B=\theta_{s_A}(A)$ such that 
$\zeta_s\vert_{B}=0$
for all $s\geq s_A$.
\end{enumerate}
Property~\ref{property_theta_becoming_stable}
implies in particular that the sets
\begin{equation*}
C_s=\{(t,q)\in [-1,\infty)\times Q
\ \vert
\ \zeta_{s'}\vert_{(t,q)}=0\ \text{for all}\ s'\geq s\}
\end{equation*}
exhaust $[-1,\infty)\times Q$,
i.e. $\bigcup_{s\in\R}\interior{C_s}=[-1,\infty)\times Q$.
This can be seen as follows. Let $B$ be
an open subset of $[-1,\infty)\times Q$
such that $\bar B$ is compact and $\interior{\bar B}=B$.
Property~\ref{property_theta_becoming_stable}
states that for $A=\theta^{-1}(\bar B)$
there is $s_A\in\R$
such that 
$\theta_s\vert_A=\theta\vert_A$
for all $s\geq s_A$
hence $\zeta_s\vert_B=0$ for all $s\geq s_A$, i.e.
$B\subset \interior{C_{s_A}}$.
Exhausting $[-1,\infty)\times Q$
by such sets $B$ yields the statement.
On the other hand, the sets
$\theta_s^{-1}(\interior{C_s})\cap ([-1,0)\times Q)$ exhaust $[-1,0)\times Q$. Observe that
$\theta\vert_{\theta_s^{-1}(\interior{C_s})}=\theta_s\vert_{\theta_s^{-1}(\interior{C_s})}$
by definition of $C_s$,
in particular $\theta^{-1}(\interior{C_s})=\theta_s^{-1}(\interior{C_s})$.
As the sets $\interior{C_s}$ exhaust
$[-1,\infty)\times Q$, the sets
$\theta^{-1}(\interior{C_s})$ exhaust 
$[-1,0)\times Q$. In summary
\begin{equation}
\begin{split}
\bigcup_{s\in\R}\interior{C_s}=[-1,\infty)\times Q\qquad\text{and}\\
\bigcup_{s\in\R}\theta_s^{-1}(\interior{C_s})\supset
 [-1,0)\times Q \ 
\end{split}\label{property_cnExhausting}
\end{equation}

\smallskip
In the second part of the proof as
the map $\varphi$ given in the theorem is assumed injective, proper and continuous
the set $[-1,\infty)\times Q$ will be seen as a subset of $M$ for simplicity. Using Equation~\eqref{equation_pointingInDirectionOfRays},  Condition~\eqref{condition_ContractionWithOmegaVanish}
requested to hold in the Theorem reads
\begin{equation*}
\omega(\zeta_s,\nu)=0\qquad \text{for all}\quad \nu\in T( (-1,\infty)\times Q)\subset TM
\end{equation*}
The next goal is to extend $\zeta_s$ to a time-dependent vector field $\xi_s$ such that $\imath_{\xi_s}\omega$ is a closed time-dependent
1-form on $M$.
In fact, it is even possible
to construct a time-dependent function
$g_s:M\to\R$ with $g_s\vert_{[-1,\infty)\times Q}=0$ such that
\begin{equation*}
\imath_{\zeta_s}\omega\vert_{(-1,\infty)\times Q}=\de g_s\vert_{(-1,\infty)\times Q}
\end{equation*}
One down to the earth way to see this might be that locally in coordinates adapted to the submanifold a local time-dependent function $g_s$ can be given explicitly. Using a partition of unity the result holds globally as the local functions $g_s$ vanish on the local submanifolds.\\
Once a neighborhood $U$ of $[0,\infty)\times Q$ is selected as stated in the Theorem, deforming $\varphi$
slightly on $[-1,0)\times Q$, it can be assumed that
$[-1,\infty)\times Q\subset U$. 
Now the time-dependent function $g_s$
can be cut off to vanish outside of $U$.
Furthermore in view of \ref{property_zeta_integrabilty},
$g_s$ can be required to satisfy
additionally that the sets
$$\bigcup_{s'\in[s-1,s+1]}\{m \in M\ \vert\ \de g_{s'}\vert_m\not=0\}$$
are relatively compact for all
$s\in\R$.\\
Starting from $g_s$ a time-dependent function $f_s$ will be defined as the limit of a sequence of 
time-dependent functions
$f_{n, s}$.
Each such time-dependent function $f_{n,s}$ and $f_s$ define a time-dependent Hamitonian vector field by
\begin{equation*}
\begin{array}{cc}
\imath_{\xi_{n,s}}\omega=\de f_{n,s}&\text{and}\\
\imath_{\xi_{s}}\omega=\de f_{s}&
\end{array}
\end{equation*}
with flows $\psi_{n,s}$ and $\psi_s$ respectively.\\
Initialize $f_{1,s}=g_s$ and set
$f_{n+1,s}=\chi_{n,s}\cdot f_{n,s}$
for a sequence of smooth time-dependent functions
$\chi_{n,s}:M\to [0,1]$
whose properties will be specified below. Since 
for all time-dependent functions $h_s\in\{f_{n,s}, f_s\}$ the set
\begin{equation*}
\bigcup_{s'\in[s-1,s+1]}\{m \in M\ \vert\ \de h_{s'}\vert_m\not=0\}\ \subset
\bigcup_{s'\in[s-1,s+1]}\{m \in M\ \vert\ \de g_{s'}\vert_m\not=0\}
\end{equation*}
is relatively compact for all
$s\in\R$, the flows $\psi_{n,s}$ and $\psi_s$ are defined for all $s\in\R$ globally.\\
In order to define $\chi_{n,s}$
let $L_n$ be an exhausting sequence of compact subsets of $M$, i.e. $L_0=\emptyset$, $L_n\subset\interior{L_{n+1}}$ and $\bigcup_{n\in\N}L_n=M$. In view of \eqref{property_cnExhausting},
this choice can be made such that $L_n\cap [-1,\infty)\times Q\subset\interior{C_n}$ for all 
$n\in\N$. Analogously using \eqref{property_cnExhausting} again, an exhaustion $K_n$ of $M\backslash ([0,\infty)\times Q)$ is chosen such that $K_n\cap [-1,0)\times Q\subset \theta_n^{-1}(\interior{C_n})$.\\
The time-dependent functions $\chi_{n,s}$ shall satisfy
\begin{enumerate}[label=($\chi$\arabic*)]
\item $\chi_{n, s}(m)=1$ if $n\geq s$,
\label{condition_chi_convergenceInN}
\item $\chi_{n,s}(m)=1$ if $m\in[-1,\infty)\times Q$ and $f_{n,s}(m)\not=0$.
\end{enumerate}
Property~\ref{condition_chi_convergenceInN}
implies that 
\begin{equation*}
\begin{split}
f_{n,s}=f_{m,s}\quad\text{and}\quad 
\xi_{n,s}=\xi_{m,s}\\\text{and hence}\quad 
\psi_{n,s}=\psi_{m,s}\quad\mbox{ }
\end{split}
\qquad\text{for all}\quad n,m\geq s\ . 
\end{equation*}
So, $f_{n,s}$ converges to $f_s$
for $n\to\infty$,
and consequently $\xi_{n,s}$ converges to $\xi_s$ and
so the sequence of corresponding flows
$\psi_{n,s}$ converges to the flow $\psi_s$ of $\xi_s$. The limits satisfy
\begin{equation*}
f_{n,s}=f_{s}\quad\text{and}\quad 
\xi_{n,s}=\xi_{s}\quad\text{and}\quad 
\psi_{n,s}=\psi_{s}
\qquad\text{for all}\quad n\geq s\ . 
\end{equation*}
As the choice of $\chi_{n,s}$ does not change $\psi_{s'}$ for all $s'\leq n$, in particular
$\psi_n$, the following third
condition on $\chi_{n,s}$ can be required,
since the sequence $\chi_{n,s}$ can be defined inductively.
\begin{enumerate}[label=($\chi$\arabic*),start=3]
\item $\chi_{n,s}(m)=0$ if $s\in [n+\delta, \infty)$ and $m\in L_n\cup \psi_n(K_n)$\label{condition_chi_ns_ImpliesBecomingStableForLargeS}
\label{condition_chi_cutoffArea}
\end{enumerate}
By making use of Property~\ref{property_zeta_constantAroundN} on the side, this condition~\ref{condition_chi_ns_ImpliesBecomingStableForLargeS} implies
\begin{equation*}
\psi_{s}\vert_{K_n}=\psi_{s'}\vert_{K_n}\quad\text{for all}\ s, s'\geq n\ .
\end{equation*}
Thus, since $K_n$ has been chosen to exhaust
$M_0=M\setminus ([0,\infty)\times Q)$,
the diffeomorphisms $\psi_s$ converge 
for $s\to\infty$
on $M_0$ to a diffeomorphism
$\psi:M_0\to\psi(M)\subset M$.
Recall that each diffeomorphism
$\psi_{n,s}$ as well as each diffeomorphism
$\psi_s$ is in fact a symplectomorphism, so the limit $\psi$ satisfies
$\psi^*\omega=\omega$.\\
Furthermore all time-dependent functions $g_s$, $f_{n,s}$ and $f_s$ and therefore all
time-dependent vector fields $\xi_{n,s}$ and $\xi_s$ are constructed to vanish
outside $U$.
Thus all flows $\psi_{n,s}$ and $\psi_s$ and a forteriori $\psi$ equal the identity outside $U$.\\
Finally $\psi$ will be shown to be surjective to $M$. By construction $\psi\vert_{[-1,0)\times Q}=\theta$ the image of $\theta$, namely $[-1,\infty)\times Q$, is contained in the image of $\psi$. For a compact set $L$ disjoint to
$[-1,\infty)\times Q$
there is $n\in\N$ such that $L\subset L_n$. Furthermore $\psi_n$ is a bijection of $M$ and the restriction
$\psi_n\vert_{M\setminus ([-1,\infty)\times Q)}$ a bijection of ${M\setminus ([-1,\infty)\times Q)}$, i.e. $L\subset \psi_n({M\setminus ([-1,\infty)\times Q)})$. For each $l\in L$ there is an
$m\in\psi_n^{-1}(L)$ such that
\begin{equation*}
l=\psi_{s}(m)=
\psi_{s'}(m)
\quad\text{for all}\ s, s'\geq n
\end{equation*}
with $m\not\in [-1,\infty)\times Q$
and hence as $\psi(m)=\psi_s(m)$ for all $s\geq n$,
the equation $\psi(m)=l$ holds
which
finishes the proof showing
surjectivity of $\psi:M_0\to M$.
\end{proof}
\paragraph{Acknowledgement}
The author expresses his gratitude
to A.~Weinstein for drawing his interest to the question
and to the result of X.~Tang (\cite{Tan18}) in the discussion process of \cite{Str20}.
Furthermore he likes to thank A.~Abbondandolo.

\bigskip
{\footnotesize
The research for this work was partially supported by the SFB/Transregio 191 “Symplektische Strukturen in Geometrie, Algebra und Dynamik”.}

\bibliographystyle{alpha}

\begin{thebibliography}{McD87}

\bibitem[Gro85]{Gro85}
M.~Gromov.
\newblock Pseudo holomorphic curves in symplectic manifolds.
\newblock {\em Invent. Math.}, 82(2):307--347, 1985.

\bibitem[McD87]{McD87}
D.~McDuff.
\newblock Symplectic structures on {${\mathbb{R}}^{2n}$}.
\newblock In {\em Aspects dynamiques et topologiques des groupes infinis de
  transformation de la m\'{e}canique ({L}yon, 1986)}, volume~25 of {\em Travaux
  en Cours}, pages 87--94. Hermann, Paris, 1987.

\bibitem[MT93]{McDTra93}
D.~McDuff and L.~Traynor.
\newblock The {$4$}-dimensional symplectic camel and related results.
\newblock In {\em Symplectic geometry}, volume 192 of {\em London Math. Soc.
  Lecture Note Ser.}, pages 169--182. Cambridge Univ. Press, Cambridge, 1993.

\bibitem[Str20]{Str20}
B.~Stratmann.
\newblock Nowhere vanishing primitive of a symplectic form.
\newblock arXiv:2003.02602, 2020.

\bibitem[Tan18]{Tan18}
X.~Tang.
\newblock Removing a ray from a noncompact symplectic manifold.
\newblock arXiv:1812.00444, 2018.

\bibitem[Tra93]{Tra93}
L.~Traynor.
\newblock Symplectic embedding trees for generalized camel spa\-ces.
\newblock {\em Duke Math. J.}, 72(3):573--594, 1993.

\end{thebibliography}

\textsc{Department of Mathematics\\
Ruhr-Universit\"at Bochum\\
44780 Bochum\\
Germany}

\smallskip
\textit{E-mail address} \texttt{bernd.x.stratmann@rub.de}
\end{document}